# On the generalized Bernoulli numbers

Donal F. Connon

dconnon@btopenworld.com

21 January 2018

**Abstract**

We derive an expression for the generalized Bernoulli numbers in terms of the Bernoulli numbers involving the (exponential) complete Bell polynomials.

## 1. Introduction

The Bernoulli numbers $B_n$ are defined by the generating function [6, p.59]

(1.1) $$\frac{x}{e^x-1} = \sum_{n=0}^{\infty} \frac{B_n}{n!} x^n$$

and the generalized Bernoulli numbers $B_n^{(\alpha)}$ are defined by the generating function [6, p.61]

(1.2) $$\left(\frac{x}{e^x-1}\right)^{\alpha} = \sum_{n=0}^{\infty} \frac{B_n^{(\alpha)}}{n!} x^n$$

and we see that $B_n^{(1)} = B_n$.

In the UK, we would spell "generalized" as generalised. However, a Google search using the American spelling yields vastly more results and I have therefore decided to employ that spelling in this paper.

## 2. The (exponential) complete Bell polynomials

The (exponential) complete Bell polynomials may be defined by $Y_0 = 1$ and for $r \geq 1$

(2.1) $$Y_r(x_1,...,x_r) = \sum_{\pi(r)} \frac{r!}{k_1! k_2!... k_r!} \left(\frac{x_1}{1!}\right)^{k_1} \left(\frac{x_2}{2!}\right)^{k_2} ... \left(\frac{x_r}{r!}\right)^{k_r}$$

where the sum is taken over all partitions $\pi(r)$ of $r$, i.e. over all sets of integers $k_j$ such that

$$k_1 + 2k_2 + 3k_3 + \cdots + rk_r = r$$

The definition (2.1) immediately implies the following relation

(2.2) $\quad Y_r(ax_1, a^2 x_2, ..., a^r x_r) = a^r Y_r(x_1,...,x_r)$

and with $a = 1$ we have

(2.3) $\qquad Y_r(-x_1, x_2, ..., (-1)^r x_r) = (-1)^r Y_r(x_1, ..., x_r)$

The complete Bell polynomials have integer coefficients and the first six are set out below (Comtet [4, p.307])

(2.4) $\qquad\qquad\qquad Y_1(x_1) = x_1$

$$Y_2(x_1, x_2) = x_1^2 + x_2$$

$$Y_3(x_1, x_2, x_3) = x_1^3 + 3x_1 x_2 + x_3$$

$$Y_4(x_1, x_2, x_3, x_4) = x_1^4 + 6x_1^2 x_2 + 4x_1 x_3 + 3x_2^2 + x_4$$

$$Y_5(x_1, x_2, x_3, x_4, x_5) = x_1^5 + 10x_1^3 x_2 + 10x_1^2 x_3 + 15x_1 x_2^2 + 5x_1 x_4 + 10x_2 x_3 + x_5$$

$$Y_6(x_1, x_2, x_3, x_4, x_5, x_6) = x_1^6 + 6x_1 x_5 + 15x_2 x_4 + 10x_2^3 + 15x_1^2 x_4 + 15x_2^3 + 60x_1 x_2 x_3$$

$$+ 20x_1^3 x_3 + 45x_1^2 x_2^2 + 15x_1^4 x_1 + x_6$$

Comtet [4, p.134] shows that the complete Bell polynomials are also given by the exponential generating function

(2.5) $\qquad\qquad \exp\left(\sum_{j=1}^{\infty} x_j \frac{t^j}{j!}\right) = \sum_{n=0}^{\infty} Y_n(x_1, ..., x_n) \frac{t^n}{n!}$

We note that

$$\sum_{n=0}^{\infty} Y_n(ax_1, ..., ax_n) \frac{t^n}{n!} = \exp\left(\sum_{j=1}^{\infty} ax_j \frac{t^j}{j!}\right)$$

$$= \exp\left[a\left(\sum_{j=1}^{\infty} x_j \frac{t^j}{j!}\right)\right]$$

$$= \left[\exp\left(\sum_{j=1}^{\infty} x_j \frac{t^j}{j!}\right)\right]^a$$

and thus we have

(2.6) $\qquad \left[\sum_{n=0}^{\infty} Y_n(x_1, ..., x_n) \frac{t^n}{n!}\right]^a = \sum_{n=0}^{\infty} Y_n(ax_1, ..., ax_n) \frac{t^n}{n!}$

Let us now consider a function $f(t)$ which has a Taylor series expansion around $x$: we have

$$e^{f(x+t)} = \exp\left(\sum_{j=0}^{\infty} f^{(j)}(x)\frac{t^j}{j!}\right) = e^{f(x)}\exp\left(\sum_{j=1}^{\infty} f^{(j)}(x)\frac{t^j}{j!}\right)$$

$$= e^{f(x)}\sum_{n=0}^{\infty} Y_n\left(f^{(1)}(x), f^{(2)}(x),..., f^{(n)}(x)\right)\frac{t^n}{n!}$$

where we have used (2.5). We see that

$$\frac{d^r}{dx^r}e^{f(x)} = \frac{\partial^r}{\partial x^r}e^{f(x+t)}\bigg|_{t=0}$$

and hence we have due to symmetry

$$\frac{d^r}{dx^r}e^{f(x)} = \frac{\partial^r}{\partial t^r}e^{f(x+t)}\bigg|_{t=0}$$

We therefore obtain

(2.7) $$\frac{d^r}{dx^r}e^{f(x)} = e^{f(x)}Y_r\left(f^{(1)}(x), f^{(2)}(x),..., f^{(r)}(x)\right)$$

as noted by Kölbig [5] and Coffey [2].

Suppose that $h'(x) = h(x)g(x)$ and let $f(x) = \log h(x)$. We see that

$$f'(x) = \frac{h'(x)}{h(x)}$$

and, since $h(x) = e^{f(x)}$, using (2.7) we have

(2.8) $$\frac{d^r}{dx^r}h(x) = \frac{d^r}{dx^r}e^{\log h(x)} = h(x)Y_r\left(g(x), g^{(1)}(x),..., g^{(r-1)}(x)\right)$$

where we define $g(x) := f'(x)$. We then have the Maclaurin expansion

(2.9) $$h(x) = h(0)\sum_{n=0}^{\infty} Y_r\left(g(0), g^{(1)}(0),..., g^{(r-1)}(0)\right)\frac{x^r}{r!}$$

**3. An expression for the Bernoulli numbers**

We see from (1.1) that

$$\frac{1}{e^t - 1} = \sum_{n=0}^{\infty} \frac{B_n}{n!} t^{n-1}$$

and integration gives us

$$\int_a^x \frac{dt}{e^t - 1} = \int_a^x \frac{e^{-t} dt}{1 - e^{-t}} = \log(1 - e^{-x}) - \log(1 - e^{-a})$$

$$= \log x - \log a + \sum_{n=1}^{\infty} \frac{B_n}{n \, n!} [x^n - a^n]$$

Hence, we obtain

$$\log(1 - e^{-x}) - \log(1 - e^{-a}) = \log x - \log a + \sum_{n=1}^{\infty} \frac{B_n}{n \, n!} [x^n - a^n]$$

We may express this as

$$\log \frac{1 - e^{-x}}{x} = \log \frac{1 - e^{-a}}{a} + \sum_{n=1}^{\infty} \frac{B_n}{n \, n!} [x^n - a^n]$$

Using L'Hôpital's rule we see that

$$\lim_{a \to 0} \log \frac{1 - e^{-a}}{a} = 0$$

and with $x \to -x$ we obtain

(3.1) $$\log\left(\frac{x}{e^x - 1}\right) = \sum_{n=1}^{\infty} \frac{(-1)^{n+1} B_n}{n \, n!} x^n$$

which was proved by Ramanujan for $|x| < 2\pi$ [1, p.119].

Let us consider the function $\log h(x)$ with the following Maclaurin expansion

(3.2) $$\log h(x) = b_0 + \sum_{n=1}^{\infty} \frac{b_n}{n} x^n$$

Then reference to (2.9) shows that

(3.3) $$h(x) = e^{b_0} \sum_{n=0}^{\infty} Y_n \left(0! b_1, 1! b_2, \ldots, (n-1)! b_n\right) \frac{x^n}{n!}$$

Hence we have using (3.1)

(3.4) $$\frac{x}{e^x-1} = \sum_{n=0}^{\infty} Y_n\left(\frac{B_1}{1}, -\frac{B_2}{2}, ..., (-1)^{n+1}\frac{B_n}{n}\right)\frac{x^n}{n!}$$

and we therefore see that

$$\sum_{n=0}^{\infty} \frac{B_n}{n!} x^n = \sum_{n=0}^{\infty} Y_n\left(\frac{B_1}{1}, -\frac{B_2}{2}, ..., (-1)^{n+1}\frac{B_n}{n}\right)\frac{x^n}{n!}$$

Equating coefficients gives us

(3.5) $$B_n = Y_n\left(\frac{B_1}{1}, -\frac{B_2}{2}, ..., (-1)^{n+1}\frac{B_n}{n}\right)$$

We have [3, p.415]

(3.6) $$Y_{r+1}(x_1, ..., x_{r+1}) = \sum_{k=0}^{r} \binom{r}{k} Y_{r-k}(x_1, ..., x_{r-k}) x_{k+1}$$

and employing (3.5) we obtain the recurrence relation

(3.7) $$B_{r+1} = \sum_{k=0}^{r} (-1)^k \binom{r}{k} \frac{B_{r-k} B_{k+1}}{k+1}$$

which may also be expressed as

(3.7) $$rB_r = \sum_{k=1}^{r} (-1)^{k+1} \binom{r}{k} B_{r-k} B_k$$

## 4. An expression for the generalized Bernoulli numbers

With regard to (3.4) we have

$$\left(\frac{x}{e^x-1}\right)^\alpha = \left[\sum_{n=0}^{\infty} Y_n\left(\frac{B_1}{1}, -\frac{B_2}{2}, ..., (-1)^{n+1}\frac{B_n}{n}\right)\frac{x^n}{n!}\right]^\alpha$$

and employing (2.6) we obtain

$$\left[\sum_{n=0}^{\infty} Y_n\left(\frac{B_1}{1}, -\frac{B_2}{2}, ..., (-1)^{n+1}\frac{B_n}{n}\right)\frac{x^n}{n!}\right]^\alpha = \sum_{n=0}^{\infty} Y_n\left(\frac{\alpha B_1}{1}, -\frac{\alpha B_2}{2}, ..., (-1)^{n+1}\frac{\alpha B_n}{n}\right)\frac{x^n}{n!}$$

Then using the definition (1.2) we get

(4.1) $$B_n^{(\alpha)} = Y_n\left(\frac{\alpha B_1}{1}, -\frac{\alpha B_2}{2}, ..., (-1)^{n+1}\frac{\alpha B_n}{n}\right)$$

Srivastava and Todorov [6, p.62] have determined a formula for $B_n^{(\alpha)}$.

**5. Open access to our own work**

This paper contains references to various other papers and, rather surprisingly, all of them are currently freely available on the internet. Surely now is the time that <u>all</u> of <u>our</u> work should be freely accessible by <u>all</u>. The mathematics community should lead the way on this by publishing <u>everything</u> on arXiv, or in an equivalent open access repository. We think it, we write it, so why hide it? You know it makes sense.

Wessex House,
Devizes Road,
Upavon,
Pewsey,
Wiltshire SN9 6DL